\documentclass[11pt,a4paper]{article}
\usepackage[latin1]{inputenc}
\usepackage{amsmath}
\usepackage{amsfonts}
\usepackage{amssymb}
\usepackage{graphicx}

\newcommand{\D}{\mathcal{D}}

\newcommand{\Sg}{\mathcal{S}}

\newcommand{\dto}{\frac{d}{dt}_{|t=0}}

\newcommand{\ad}{\mbox{ad}}

\newcommand{\demo}{\noindent\textit{Proof. }}

\newcommand{\A}{\mathcal{A}}
\newcommand{\Z}{\mathcal{Z}}
\newcommand{\N}{\mathbb{N}}

\newcommand{\R}{\mathbb{R}}
\newcommand{\h}{\mathfrak{h}}
\newcommand{\hn}{\mathfrak{h}^n}
\newcommand{\Hh}{\mathbb{H}}
\newcommand{\Hn}{\mathbb{H}^n}
\newcommand{\g}{\mathfrak{g}}
\newcommand{\LL}{\mathcal{L}}
\newcommand{\LSS}{\mathcal{LSS}}
\newcommand{\LS}{\mathcal{LS}}
\newcommand{\Li}{\mathcal{L}}
\newcommand{\Lo}{\mathcal{L}_0}
\newcommand{\X}{\mathcal{X}}

\newcommand{\ds}{\displaystyle}

\begin{document}

\title{ Controllability of Linear Systems on Generalized Heisenberg Groups}
\author{Mouhamadou Dath\footnote{Universit\'e de Dakar Senegal   E-mail: rassouldath@yahoo.fr }\ \ {\small and} 
 Philippe Jouan\footnote{Lab. R.~Salem, CNRS UMR 6085, Universit\'e
    de Rouen, avenue de l'universit\'e BP 12, 76801
    Saint-\'Etienne-du-Rouvray France. E-mail: Philippe.Jouan@univ-rouen.fr}}
\maketitle
\begin{abstract}
This paper is devoted to the study of controllability of linear systems on generalized Heisenberg groups.

Some general necessary controllability conditions and some sufficient ones are provided.

We introduce the notion of decoupled systems, and more precise controllability criteria are stated for them.

\vskip 0.2cm


\vskip 0.2cm
\end{abstract}

\section{Introduction}\label{introd}

In this paper we are interested in controllability properties of linear systems on generalized Heisenberg groups, which are controlled systems
$$
(\Sigma) \qquad \dot{g}=\X(g)+\sum_{j=1}^m u_jB_j(g)
$$
where $\X$ is a linear vector field, that is a vector field whose flow is a one-parameter group of automorphisms, and the $B_j$'s are right-invariant.

This paper follows and generalizes \cite{DJ14} where necessary and sufficient controllability conditions on the $2$-dimensional affine group and the $3$-dimensional Heisenberg group were stated.

As well as in \cite{DJ14} we are interested in controllability for unbounded inputs. Controllability for bounded inputs is also an issue of interest, it has been studied by Adriano da Silva in \cite{Dasilva14} and \cite{Dasilva15}.

Since the Heisenberg groups are nilpotent we should recall the following theorem: (\cite{Jouan11}) \textit{if the derivation associated to the linear field is inner (see Section \ref{vectfieldprelim} for the definition) then System $(\Sigma)$ is controllable if and only if the Lie algebra generated by the controllable vector fields $B_1,\dots,B_m$ is equal to the Lie algebra of $G$}. This theorem has the consequence that on a nilpotent Lie group the derivation associated to a linear system of interest is not inner, and the system cannot be associated to an invariant one as in the semi-simple case (see \cite{Jouan11})). Moreover the methods used for proving controllability or non controllability are different from the ones used on semi-simple Lie groups (see \cite{Jouan16}).

A large part of our proofs makes use of subgroups and quotients. The subgroup in consideration is often the derived one but we are more generally interested in closed subgroups that are invariant under the flow of the linear vector field. The linear system can then be projected to the quotient group, and our approach is to rely the controllability properties of the system to the ones of the systems induced on the involved subgroup and on the related quotient group.

\vskip 0.2cm

The paper is organized as follows. In Section \ref{Basic} we first recall basic definitions and facts concerning linear systems on Lie groups. Then the relations between the controllability properties of the system and the ones of the systems induced on some subgroup and on the related quotient group are recalled and the section finishes by the 
statement of the necessary and sufficient controllability conditions on the $3$-dimensional Heisenberg group.

In Section \ref{LinSysHn} the derivations of the Lie algebra $\hn$ and  their expressions in particular basis that we call symplectic are analyzed. The associated linear vector fields are computed.

Some general necessary controllability conditions and some sufficient ones are stated in Section \ref{NecAndSufCond}.

To finish we consider what we call decoupled systems in Section \ref{Decoupling}. These systems are easier to analyze than the completely general ones, and we obtain more accurate results, in particular necessary and sufficient controllability conditions for decoupled systems in $\Hh^2$.


\section{Basic definitions and known results}\label{Basic}

More details about linear vector fields and linear systems can be found in \cite{Jouan09} and \cite{Jouan11}.

\subsection{Linear vector fields}\label{vectfieldprelim}

Let $G$ be a connected Lie group and $\g$ its Lie algebra (the set of right-invariant vector fields, identified with the tangent space at the identity). A vector field on $G$ is said to be {\it linear} if its flow is a one-parameter group of automorphisms. Actually the linear vector fields are nothing else than the so-called infinitesimal automorphisms in the Lie group litterature (see \cite{Bourbaki2} for instance). They can also be characterized as follows:

{\it A vector field $\X$ on a connected Lie group $G$ is linear if and only if 
it belongs to the normalizer of $\g$ in the algebra of analytic vector fields of $G$,
that is
$$
\forall Y\in \g \qquad [\X,Y]\in\g,
$$
and verifies $\X(e)=0$.}

On account of this characterization, one can associate to a linear vector field $\X$ the derivation $D=-ad(\X)$ of the Lie algebra $\g$ of $G$.
The minus sign in this definition comes from the formula
$[Ax,b]=-Ab$ in $\R^n$. It also enables to avoid a minus sign in the useful formula:
\begin{equation}\label{equationfonda}
\forall Y\in\g, \quad \forall t\in \R \qquad \varphi_t(\exp
Y)=\exp(e^{tD}Y).
\end{equation}

In the case where this derivation is inner, that is $D=-ad(X)$ for some right-invariant vector field $X$ on $G$, the linear vector field splits into $\X=X+\mathcal{I}_*X$, where $\mathcal{I}$ stands for the diffeomorphism $g\in G\longmapsto \mathcal{I}(g)=g^{-1}$. Thus $\X$ is the sum of the right-invariant vector field $X$ and the left-invariant one $\mathcal{I}_*X$.

About the existence of linear vector fields, we have:
\newtheorem{champ linearsimplyconnex}{Theorem}
\begin{champ linearsimplyconnex}\label{champ linearsimplyconnex}
\textnormal{(\cite{Jouan09})} \textit{The group G is assumed to be (connected and) simply connected. Let D be a derivation of its
Lie algebra g. Then there exists one and only one linear vector field on G whose associated derivation is D.}
\end{champ linearsimplyconnex}

Throughout the paper the flow of a linear vector field $\X$ will be denoted by $(\varphi_t)_{t\in\R}$.


\subsection{Linear systems}\label{linearsystem}
\newtheorem{syslin}{Definition}
\begin{syslin}
A {\bf linear system} on a connected Lie group $G$ is a controlled system
$$(\Sigma) \qquad \dot{g}=\X(g)+\sum_{j=1}^m u_jB_j(g)$$
where $\X$ is a linear vector field and the $B_j$'s are right-invariant ones. The control $u=(u_1,\dots,u_m)$ takes its values in $\R^m$.
\end{syslin}

An input $u$ being given (measurable and locally bounded), the corresponding trajectory of $(\Sigma)$ starting from the identity $e$ will be denoted by $e_u(t)$, and the one starting from the point $g$ by $g_u(t)$. A straightforward computation shows that
$$
g_u(t)=e_u(t)\varphi_t(g).
$$

\vskip 0.2cm

\noindent{\bf Notations}.
We denote by $\A(g,t)=\{g_u(t);\ u\in L^{\infty}[0,t]\}$ (resp. $\A(g,\leq t)$)  (resp. $\A(g)$) the reachable set from $g$ in time $t$ (resp. in time less than or equal to $t$) 
(resp. in any time). In particular the reachable sets from the identity $e$ are denoted by
$$
\A_t=\A(e,t)=\A(e,\leq t) \qquad\mbox{ and }\qquad \A=\A(e).
$$

We also denote by  $ \A^{-}=\{ g \in G; \  e \in \A(g)\}$ the set of points from which the identity can be reached. It is equal to the attainability set from the identity for the time-reversed system.


\subsubsection{The system Lie algebra and the rank condition}

Let $V$ stand for the subspace of $\g$ generated by $\{B_1,\dots,B_m\}$, let us denote by $DV$ the smallest $D$-invariant subspace of $\g$ that contains $V$, i.e. $DV=\mbox{Span}\{D^kY;\ Y\in V \mbox{ and } k\in\N\}$,
and let $\mathcal{LA}(DV)$ be the $\g$ subalgebra generated by $DV$ (as previously $D=-\ad(\X)$).

\newtheorem{AlgebreLie}{Proposition}
\begin{AlgebreLie}\label{AlgebreLie}
The subalgebra $\mathcal{LA}(DV)$ of $\g$ is $D$-invariant. It is therefore equal to the zero-time ideal $\Lo$, and the system Lie algebra $\Li$ is equal to
$$
\R\X\oplus \mathcal{LA}(DV)=\R\X\oplus \Lo.
$$
The rank condition is satisfied by $(\Sigma)$ if and only if $\Lo=\g$.
\end{AlgebreLie}


\subsubsection{The Lie saturate}\label{Sature}

The Lie saturate $\LS(\Sigma)$ of $(\Sigma)$ (resp. the strong Lie saturate $\LSS(\Sigma)$ of $(\Sigma)$) is the set of vector fields $f$ belonging to the system Lie algebra $\Li$ and whose flow $(\phi_t)_{t\in\R}$ satisfies
$$
\forall g\in G,\ \forall t\geq 0\qquad \phi_t(g)\in\overline{\A(g)} \qquad (\mbox{resp. } \phi_t(g)\in\overline{\A(g,\leq t)})
$$
as soon as $\phi_t(g)$ is defined (see \cite{Jurdjevic97}). Here $\overline{\A(g)}$ stands for the closure of $\A(g)$.

As a first consequence of the enlargement technics, that consists to add to a system some vector fields of the (strong) Lie saturate, we have the following proposition:

\newtheorem{Extension1}[AlgebreLie]{Proposition}
\begin{Extension1}\label{Extension1}
Let $\h$ be the subalgebra of $\g$ generated by $\{B_1,\dots,B_m\}$. It is included in $\LSS(\Sigma)$, so that $(\Sigma)$ can be enlarged to the system
$$
(\widetilde{\Sigma}) \qquad \dot{g}=\X(g)+\sum_{j=1}^p u_j\widetilde{B}_j(g),
$$
where $\widetilde{B}_1,\dots,\widetilde{B}_p$ is a basis of $\h$,
without modifying the sets
$
\overline{\A(g,\leq t)}.
$
\end{Extension1}


\subsubsection{Local controllability and the ad-rank condition}\label{ContLocCondAlg}

It is well known that a system is locally controllable at an equilibrium point as soon as the linearized system is controllable (see \cite{NvdS90} for instance). 
In this assertion "locally controllable" at a point $g$ means that the set $\A(g,t)$ is a neighbourhood of $g$ for all $t>0$.

According to Proposition \ref{Extension1}, we can consider the linearization of the extended system $(\widetilde{\Sigma})$ instead of the one of the system itself. This leads to the following definition:

\newtheorem{Rangalgebrique}[syslin]{Definition}
\begin{Rangalgebrique}
System $(\Sigma)$ is said to satisfy the {\bf ad-rank condition} if $D\h=\g$, in other words if the linearized system of $(\widetilde{\Sigma})$ is controllable.
\end{Rangalgebrique}
As the rank condition is implied by the ad-rank one we obtain at once:
\newtheorem{Contlocal}[AlgebreLie]{Proposition}
\begin{Contlocal}
If the ad-rank condition is satisfied then for all $t>0$ the reachable set $\A_t$ is a neighbourhood of $e$.\\
\end{Contlocal}


\subsubsection{Equivalence of systems}\label{EquivlSys}

The next proposition \ref{theohei2} comes from \cite{DJ14}. It shows that for a control system, the fact to be linear is preserved by Lie group automorphims. Actually we will see in Proposition \ref{Auto.Sympl} (see Section \ref{SymplecticBasis}), that transforming a linear system by automorphism is equivalent to choosing a suitable basis of a particular type.

\newtheorem{theohei2} [AlgebreLie]{Proposition}
\begin{theohei2}\label{theohei2}
Let $(\Sigma)$: $\dot{g}=\X(g)+\sum_{j=1}^m u_jB_j(g)$ be a linear system on a connected and simply connected Lie group $G$, $D$ the derivation associated to $\X$ and $ P $  an  automorphism of $ \g $. Then $ (\Sigma) $ is equivalent by group automorphism to
$$
(\widetilde{\Sigma})\qquad \dot{g}=\widetilde{\X _g}+\sum_{j=1}^m u_jPB_j(g)
$$ 
where $\widetilde{\X}$ is the linear vector field whose associated  derivation is $\widetilde{D} = PDP^{-1}$.
\end{theohei2}


\subsection{Controllability and quotient groups}

It is a well known fact that the derived subalgebras, as well as the subalgebras of the lower central series, are characteristic ideals of $\g$, hence invariant by derivations. The corresponding connected subgroups are therefore normal and invariant by the flow of $\X$ (see \cite{Jouan09}). Moreover these subgroups are closed when $G$ is simply connected (see \cite{Bourbaki2}). Notice also that in the case when $G$ is simply connected, the group $G/\D^1G$ is abelian and simply connected, hence diffeomorphic to $\R^n$ for some $n$, and that the induced system on $G/\D^1G$ is linear in the classical sense. 

The two forthcoming propositions are proved in \cite{DJ14}.

\newtheorem{propgen1}[AlgebreLie]{Proposition}
\begin{propgen1}\label{propgen1}
Let $H$ be a closed  subgroup of $G$, globally invariant under the flow of $\X$.

The linear system $(\Sigma)$ on $G$, assumed to satisfy the rank condition, is controllable if and only if both conditions hold:
\begin{enumerate}
	\item the system induced on $G/H$ is controllable;
	\item the subgroup $H$ is included in the closures  $\overline{\A}$ and $\overline{\A^{-}}$  of $\A$ and $\A^{-}$.
\end{enumerate}
\end{propgen1}

\newtheorem{propgen2}[AlgebreLie]{Proposition}
\begin{propgen2}\label{propgen2}
Let us assume that $ H $ is a connected and closed subgroup of $G$ and that the restriction of $ \X $ to $ H $  vanishes. Then $ H $ is included in $\A$ (resp. in $ \A^{-} $) if and only if $ \A \cap H $(resp. $\A^{-}\cap H $) is a neighbourhood of $ e $ in $ H $.
\end{propgen2}

\vskip 0.2cm

\noindent \textbf{Singular and regular systems.} The linear systems for which $\X$ vanishes on a connected, closed and $\X$-invariant (non trivial) subgroup will be referred to as \textit{singular systems}, the other ones being \textit{regular}. The previous proposition \ref{propgen2} is crucial in the singular case.


\subsection{Controllability on $\Hh^1$}\label{ContHun}

We consider here a system $(\Sigma)$ with one input in $\Hh^1$. It is proved in \cite{DJ14} that under the rank condition it is equivalent by automorphism to a system in "normal form" in the canonical basis, that is:

$$
(\Sigma)\qquad\dot{g}=\X(g)+ u X(g)
\qquad\mbox{where}\qquad D=\begin{pmatrix}0&b&0\\1&d&0\\0&f&d\end{pmatrix}
$$
is the matrix in the basis $(X,Y,Z)$ of the derivation associated to $\X$. Notice that the controlled vector field is the first element of that basis.

The main result of \cite{DJ14}, is the following.
\newtheorem{theohei3}[champ linearsimplyconnex]{Theorem}
\begin{theohei3}\label{theohei3}
A system in normal form is controllable if and only if one of the conditions
\begin{enumerate}
	\item[(i)] $\displaystyle b < -\frac{d^{2}}{4}$,
	\item[(ii)] $d=0$ and $f \neq 0$,
\end{enumerate}
holds.
  \end{theohei3}
  
Let us denote by $(L)$ the classical linear system
$$
(L)=\quad \left\{\begin{array}{ll}
\dot{x}=& by + u\\
\dot{y}=& x +dy
\end{array} \right.
$$
induced on the quotient $G/\Z(G)$. The eigenvalues of the matrix $\begin{pmatrix}0&b\\1&d\end{pmatrix}$ are real if and only if $\ds b\geq -\frac{d^2}{4}$. On the other hand $(\Sigma)$ is singular if and only if $d=0$ and the algebraic rank condition is satisfied if and only if $f\neq 0$. Consequently Theorem \ref{theohei3} can be restated as:

\newtheorem{coroheis}[champ linearsimplyconnex]{Theorem}
\begin{coroheis}\label{coroheis}
The one-input system $(\Sigma)$ on the Heisenberg group is controllable if and only if it satisfies the rank condition and
\begin{enumerate}
	\item[(i)] in the regular case: the eigenvalues of $(L)$ are not real;
	\item[(ii)] in the singular case: the eigenvalues of $(L)$ are not real or the ad-rank condition is satisfied.
\end{enumerate}
\end{coroheis}


\section{Linear Systems on $\Hn$}\label{LinSysHn}

The generalized Heisenberg group $\mathbb{H}^{n}$ can be defined as the following $(2n+1)$-dimensional matrix subgroup of $GL(n+2, \mathbb{R})$:
$$
\mathbb{H}^{n}=\left\{\begin{pmatrix}1&y_{1}&y_{2}&...&y_{n}&z\\0&1&0&...&0&x_{1}\\0&0&1&....&0&x_{2}\\
.&.&.&.&.&.&\\.&.&.&.&.&.&\\0&0&0&...&1&x_{n}&\\
0&0&0&...&0&1\\\end{pmatrix};\ \  x_{i},y_{i},z\in\mathbb{R} \right\}.$$

Its Lie algebra $\mathfrak{h}^{n}$ is generated by the $2n+1$ right-invariant vector fields  $X_{1},\dots,X_{n}$, $Y_{1},\dots,Y_{n}$, and $Z$ defined by
$$\left\{\begin{array}{ll}
 X_{i}& = E_{i+1,n+2} \qquad 1 \leq i \leq n \\
 Y_{i}& = E_{1,i+1} + x_{i}E_{1,n+2} \qquad 1 \leq i \leq n  \\
 Z& = E_{1,n+2}
\end{array}
\right.$$
where $E_{ij}$ is the matrix whose all entries vanish, excepted the rank $i$ and column $j$ entry which is equal to $1$.

In canonical coordinates these vector fields write:    
\begin{equation}\label{canbasis}
X_{i}= {\frac{\partial}{\partial x_{i}}},\qquad Y_{i}= {\frac{\partial}{\partial y_{i}}} + x_{i} {\frac{\partial}{\partial z}},\qquad Z= {\frac{\partial}{\partial z}}.
\end{equation}


\subsection{Symplectic basis}\label{SymplecticBasis}

The basis of $\hn$ defined above satisfies the following Lie bracket relations:
\begin{equation}
[X_{i},Y_{i}]= Y_{i}X_{i}- X_{i}Y_{i} = Z \quad\mbox{ for }\quad i= 1,\dots,n
\end{equation}
and all the other brackets vanish. In particular the center of $\hn$ is generated by the field $Z$, and is equal to the derived algebra $\D^1\hn$.
We are interested in the basis of $\hn$ that satisfy the same kind of relations and in conditions under which a given set of elements of $\hn$ is part of such a basis. For this purpose we introduce what follows.

\textit{For all $X$ and $Y$ in $\hn$ we define $l(X,Y)$ as the unique real number such that $[X,Y]=l(X,Y)Z$. The mapping $l$ is clearly a skew-symmetric bilinear form on $\hn$. Moreover its kernel (the set of $N\in\hn$ such that $\forall X\in\hn$ $l(X,N)=0$) is equal to the center $\D^1\hn=\R Z$ of $\hn$, so that it induces a non degenerated skew-symmetric form, that is a symplectic form, on the quotient $\hn/\D^1\hn$ which is a $2n$-dimensional vector space.}

\vskip 0.2cm

\noindent \textbf{Important remark} Up to a non-zero constant the symplectic form $l$ depends on the choice of the generator $Z$ of $\D^1\hn$, and we will sometimes have to multiply $Z$ it by some constant,  thus modiying $l$. In what follows $l$ will always be related to the $Z\in \D^1\hn$ under consideration. 

With this in mind, and with a clear abuse of language, we define what we call symplectic basis on $\hn$.

\newtheorem{Sympl.Basis}[syslin]{Definition}
\begin{Sympl.Basis}\label{Sympl.Basis}
A basis $(X_1,Y_1,\dots,X_n,Y_n,Z)$ of $\h^n$ is said to be symplectic if $[X_i,Y_i]=Z$, $i=1,\dots,n$, and  the other brackets vanish, in other words if and only if the set of  projections of $(X_1,Y_1,\dots,X_n,Y_n)$ on the quotient $\hn/\D^1\hn$ is a symplectic basis of  $\h^n/ \mathbb{R} Z=\h^n/\D^1\h^n$  w.r.t. the symplectic form $l$. 
\end{Sympl.Basis}

By a standard application of the linear algebra methods we can therefore state the two following propositions that will be very useful in the sequel:
\newtheorem{Symplec Vectors}[AlgebreLie]{Proposition}
\begin{Symplec Vectors}\label{Symplec Vectors}
\begin{enumerate}
\item
If $X\in \h^n\setminus D^{1}\h^{n}$ then there exists a symplectic basis\\
$(X_1,Y_1,\dots, X_n ,Y_n, Z)$ of $\hn$ such that $X = X_1$.
\item
If $X, Y \in \h^n$ and $ \left[X, Y\right] \neq 0 $ then there exists a symplectic basis
$(X_1,Y_1,\dots, X_n ,Y_n, Z)$ of $\hn$ such that $X = X_1$ and $Y = Y_1$.
\item
If $B_1,\dots,B_m$ are linearly independant in $\h^n/\D^1\h^n$ and $\left[ B_i, B_j\right] = 0$ for all $i, j = 1,\cdots,m $, then there exists a symplectic basis
$(X_1,Y_1,\dots, X_n ,Y_n, Z)$ of $\hn$ such that  $B_i = X_i$ for $i =1,\cdots, m$.  
\end{enumerate}
\end{Symplec Vectors}

The next proposition is obvious, its proof is omitted.
\newtheorem{Auto.Sympl}[AlgebreLie]{Proposition}
\begin{Auto.Sympl}  \label{Auto.Sympl}
A linear isomorphism $P$ of $\hn$ is a Lie algebra automorphism if and only if the image by $P$ of any symplectic basis of $\hn$ is again a symplectic basis of $\hn$.
\end{Auto.Sympl}

Thanks to this proposition it is equivalent to write a given system in a suitable symplectic basis or to transform it by automorphism to an equivalent system in the canonical basis. Morover we can associate to any symplectic basis a coordinate system in which the vector fields of the basis write like in (\ref{canbasis}). We make a constant use of these coordinates in the sequel.


\subsection{Derivations and linear fields on $\hn$}

\newtheorem{Derivation}[AlgebreLie]{Proposition}
\begin{Derivation}\label{Derivation}
An endomorphism $D$ of $\h^n$ is a derivation if and only if its matrix in any symplectic basis has the following form:  
$$ D= \begin{pmatrix}A_{11}&-\widetilde{A}_{21}^T& -\widetilde{A}_{31}^T &...&-\widetilde{A}_{n1}^T&0\\A_{21}&A_{22}&-\widetilde{A}_{32}^T&...&-\widetilde{A}_{n2}^T&0\\A_{31}&A_{32}&A_{33}&....&-\widetilde{A}_{n3}^T&0\\
.&.&.&.&.&.&\\
A_{n1}&A_{n2}&A_{n3}&...&A_{nn}&0\\a_{2n+1,1}&a_{2n+1,2}&a_{2n+1,3}&...&a_{2n+1,2n}&d\\
\end{pmatrix} $$
where the $A_{ij}$'s are $2\times 2$ matrices and:
\begin{itemize}
\item $\widetilde{A}_{ij}^T$ stands for the transpose of the comatrix of $A_{ij}$;
\item $tr(A_{ii}) = d$  for $1 \leq i \leq n $.
\end{itemize} 
\end{Derivation}
\demo

Let $(X_{1},Y_{1},X_{2},Y_{2}, \cdots ,X_{n},Y_{n},Z)$ be any symplectic basis $\mathcal{B}$.

First of all the one dimensional derived algebra $\D^1\hn$ being stable by the derivation $D$ we get $DZ = dZ$ for some real number $d$ (alternatively $DZ = [DX_{1}, Y_{1}]+ [X_{1}, DY_{1}] =(l(DX_1,Y_1)+l(X_1,DY_1))Z$). Consequently the matrix of $D$ writes in $\mathcal{B}$:
  $$
  \begin{pmatrix}a_{1,1}&b_{1,2}& \cdots &a_{1,2n-1}&b_{1,2n}&0\\
 a_{2,1}&b_{2,2}&\cdots &a_{2,2n-1}&b_{2,2n}&0\\
 a_{3,1}&b_{3,2}&\cdots &a_{3,2n-1}&b_{3,2n}&0\\
 a_{4,1}&b_{4,2}&\cdots &a_{4,2n-1}&b_{4,2n}&0\\
.&.&\cdots &.&.&.&\\.&.&\cdots&.&.&.&\\.&.&\cdots &.&.&.&\\
a_{2n+1,1}&b_{2n+1,2}& \cdots &a_{2n+1,2n-1}&b_{2n+1,2n}&d\\
\end{pmatrix}.
$$

In other words we have for $i=1,\dots,n$:
$$
\begin{array}{ll} DX_{i}&=\sum_{k=1}^{n}(a_{2k-1,2i-1})X_{k}+(a_{2k,2i-1})Y_{k}+(a_{2n+1,2i-1})Z\\
DY_{i}&=\sum_{k=1}^{n}(b_{2k-1,2i})X_{k}
+(b_{2k,2i})Y_{k}+(b_{2n+1,2i})Z
\end{array}
$$
\vskip 0.25cm
The following equalities hold for $i\neq\,j$: 
$$
\begin{array}{lll}
 & D[X_{i}, Y_{j}]&=[DX_{i}, Y_{j}]+[X_{i}, DY_{i}]=0  \\ 
 & D[Y_{i}, Y_{j}]&=[DY_{i}, Y_{j}]+[Y_{i}, DY_{j}]=0\\ 
 & D[X_{i}, X_{j}]&=[DX_{i}, X_{j}]+[X_{i}, DX_{j}]=0;\\ 
\mbox{and for i = j}: &&\\
 & \quad DZ&=[DX_{i}, Y_{i}] + [X_{i}, DY_{i}] =dZ. 
\end{array}
$$
From these equalities, and because the basis is symplectic, we obtain \begin{align*}
a_{2j-1,2i-1}&=-b_{2i,2j}\\
b_{2j-1,2i}&=b_{2i-1,2j}\\
a_{2i,2j-1}&=a_{2j,2i-1}\\
a_{2i-1,2i-1}+b_{2i,2i}&=d,
\end{align*} 
which finishes the proof.

\hfill $\Box$

\vskip 0.2cm

The next step consists in computing the associated linear vector field. For that purpose we introduce some notations.

\noindent\textbf{Notations}
The Lie algebra of $\mathbb{H}^{n}$ is: 
$$ \mathfrak{h}^{n} = \left\{ A=\begin{pmatrix}0&y_{1}&y_{2}&...&y_{n}&z 
\\0&0&0&...&0&x_{1}\\0&0&0&....&0&x_{2}\\
.&.&.&.&.&.&\\.&.&.&.&.&.&\\0&0&0&...&0&x_{n}&\\
0&0&0&...&0&0\\\end{pmatrix}/ x _{i}, y _{i}, z \in \mathbb{R}      \right\}$$

In what follows the elements of $\hn$ will be denoted by $A = \begin{pmatrix}
0&y_A&z_A&\\0&0&x_A&\\0&0&0&\\
\end{pmatrix}$ where $x_A=(x_1,\dots,x_n)$ and $y_A=(y_1,\dots,y_n)$ belong to $\R^n$.

The canonical scalar product of $\R^n$ will be denoted by $\langle , \rangle$, so that if $B$ is another element of $\hn$ then the matricial product $AB$ writes merely: $AB= \langle y_A, x_B \rangle Z$.

Notice also that the elements of $\Hn$ have the form $g=I+G$, where $I$ is the identity matrix of size $n+2$ and $G$ belongs to $\hn$ (but  $I+G$ is not equal to $\exp(G)$).

\newtheorem{P15}[AlgebreLie]{Proposition}
\begin{P15}\label{P15}
Let $D$ be a derivation of $\hn$. It is associated to a unique linear vector field $\X$ of $\Hn$, and $\X$ is equal, at the point $g=I+G$ of $\Hn$ to:
$$
\begin{array}{ll}
\X(g) & = DG - \frac{1}{2}dG^{2} + \frac{1}{2}(G(DG) + (DG)G)\\
       & = DG + \frac{1}{2}(\langle y_G, x_{DG}\rangle + \langle x_G , y_{DG}\rangle -d\langle x_G, y_G\rangle )Z.
\end{array}
$$
where $d$ is defined by $DZ=dZ$.
\end{P15}
\demo\\
For all $A$ in the Lie algebra $\hn$ of $\Hn$ we have $A^{2} = \langle x_A, y_A \rangle Z$ and $A^{k} = 0$ for $k>2$. Consequently the exponential mapping of $\Hn$ is:
\begin{align*}
 \exp(A)& = I_{n+2} + A + \frac{1}{2}A^{2}\\
         & = I_{n+2} + A + \frac{1}{2}\langle x_A, y_A \rangle Z,
  \end{align*}
  and its differential at the point $A$ is given by \begin{equation}\label{E1}         
T_{A} \exp.H = H + \frac{1}{2}(AH + HA).
  \end{equation}
Since the exponential map is a diffeomorphism from $\hn$ onto $\Hn$ we get for $g=I + G \in \mathbb{H}^{n}$:
\begin{equation}\label{E2}
\log(g) = G - \frac{1}{2}G^{2}.
\end{equation}
On the other hand the derivation at $t=0$ of the formula $\varphi_{t}( \exp Y) = \exp (e^{tD}Y)$ (Formula (\ref{equationfonda}) in Section \ref{vectfieldprelim}) gives:
$$
\X(\exp Y) = \dto \exp (e^{tD}Y) = T_{Y} \exp .DY,
$$
which becomes for $\exp Y = g$:
\begin{equation}\label{E3}
\X(g) = T_{\log (g)}\exp .D \log(g).
\end{equation}
Applying  \eqref{E1} in \eqref{E3}, we obtain \begin{equation}\label{E4}
\X(g) = D\log(g) + \frac{1}{2}(\log(g)(D\log(g)) +( D\log(g))\log(g)).
\end{equation}
Since $G^2=\langle x_G, y_G \rangle Z$ we have $DG^2=dG^2$ and $G^2A=AG^2=0$ for all $A\in\hn$. Consequently:
\begin{align}
D\log(g)& = DG -\frac{1}{2}DG^{2} \label{E5}=DG - \frac{1}{2}dG^{2} \\
\log(g)(D\log(g))&=(G - \frac{1}{2}G^{2})(DG - \frac{1}{2}dG^{2})=G(DG)\label{E7}\\
(D\log(g))\log(g)&= (DG - \frac{1}{2}dG^{2})(G - \frac{1}{2}G^{2})\label{E9} = (DG)G
\end{align}
Finally
$$
\begin{array}{ll}
\X(g) & = DG - \frac{1}{2}dG^{2} + \frac{1}{2}(G(DG) + (DG)G)\\
       & = DG + \frac{1}{2}(\langle y_G, x_{DG}\rangle + \langle x_G , y_{DG}\rangle -d\langle x_G, y_G\rangle )Z.
\end{array}
$$

\hfill $\Box$

\vskip 0.2cm

\textbf{An example in $\Hh^2$}. Let 
$$ D = \begin{pmatrix}
a_{11}&b_{12}&-b_{42}&b_{32}&0\\
a_{21}&d-a_{11}&a_{41}&-a_{31}&0\\
a_{31}&b_{32}&a_{33}&b_{43}&0\\
a_{41}&b_{42}&a_{43}&d-a_{33}&0\\
a_{51}&b_{52}&a_{53}&b_{54}&d\\
\end{pmatrix} \quad\mbox{ and }\ \  G = \begin{pmatrix}
x_{1}\\y_{1}\\x_{2}\\y_{2}\\z
\end{pmatrix} 
$$
Then
$$\left.\begin{array}{ll}
\X_{x_{1}}&=a_{11}x_{1}+b_{12}y_{1}-b_{42}x_{2}+b_{32}y_{2}\\
\X_{y_{1}}&=a_{21}x_{1}+(d-a_{11})y_{1}+a_{41}x_{2}-a_{31}y_{2}\\
\X_{x_{2}}&=a_{31}x_{1}+b_{32}y_{1}+a_{33}x_{2}+b_{34}y_{2}\\
\X_{y_{2}}&=a_{41}x_{1}+b_{42}y_{1}+a_{43}x_{2}+(d-a_{33})y_{2}\\
\X_{z} &=a_{51}x_{1}+b_{52}y_{1}+dz+a_{53}x_{2}+b_{54}y_{2}\\
 &+\frac{1}{2}(a_{21}x_{1}^{2}+b_{12}y_{1}^{2}+a_{43}x_{2}^{2}+b_{34}y_{2}^{2})+a_{41}x_{1}x_{2}+b_{32}y_{1}y_{2}
 \end{array}\right. $$


\section{Some necessary and some sufficient controllability conditions}\label{NecAndSufCond}

In this section we deal with the system
$$
(\Sigma) \qquad \dot{g}=\X_g+\sum_{j=1}^m u_j B_j(g)
$$
where $\X$ is a linear vector field on the $(2n+1)$-dimensional Heisenberg group $\Hn$ and the $B_j$'s are right invariant.


\subsection{Obstruction to Controllability}

It is proved in the next subsection that, among other conditions, the system is controllable as soon as the generator $Z$ of $\D^1\hn$ belongs to $\mbox{Span}\{B_1,\dots,B_m\}$. The purpose being herein to state conditions of non controllability we assume that:
$$
Z\neq \mbox{Span}\{B_1,\dots,B_m\}
$$

Thanks to that assumption, we can choose a symplectic basis $(X_1,Y_1,\dots, X_n,Y_n,Z)$ such that the $B_j$'s belong to the subspace of $\hn$ generated by the $X_i$'s and the $Y_i$'s. In the associated coordinates, the differential equation satisfied by the last coordinate $z$ does not depend on the controls, it is
$$
\dot{z}=dz+l(x_1,y_1,\dots,x_n,y_n)+Q(x_1,y_1,\dots,x_n,y_n),
$$
where $dz+l(x_1,y_1,\dots,x_n,y_n)$ is the linear form that comes from the last line of $D$ and  $Q$ is a quadratic form in the $2n$ variables $x_1,y_1,\dots,x_n,y_n$.

\newtheorem{Formequadratique}[champ linearsimplyconnex]{Theorem}
\begin{Formequadratique}\label{Formequadratique}
It is assumed that $Z\notin \LL\A\{B_1,\dots,B_m\}$ and $d\neq 0$.

If the quadratic form $Q$ is non negative, and if $\ker(l)\subset \ker(Q)$, then the system is not controllable.
\end{Formequadratique}

\demo
The assumptions on $l$ and $Q$ imply that the polynomial $l(x_1,y_1,\dots,x_n,y_n)+Q(x_1,y_1,\dots,x_n,y_n)$ has a  minimum $\mu\in\R$. Hence:
$$
\dot{z}<0 \Longleftrightarrow dz<-l(x_1,y_1,\dots,x_n,y_n)-Q(x_1,y_1,\dots,x_n,y_n)\Longrightarrow dz\leq -\mu
$$
Consequently:
$$
\begin{array}{ll}
\mbox{If } d>0, & \mbox{ then } z\geq -\frac{\mu}{d}\Longrightarrow \dot{z}\geq 0 \qquad \forall x_i,y_i\\
\mbox{If } d<0, & \mbox{ then } z\leq -\frac{\mu}{d}\Longrightarrow \dot{z}\geq 0 \qquad \forall x_i,y_i
\end{array}
$$
This proves that $(\Sigma)$ is not controllable since the hyperplane $\{\ds z=-\frac{\mu}{d}\}$ can be crossed in only one direction.

\hfill $\Box$

\vskip 0.2cm

We can actually go further by considering some particular modification of the variable $z$. \textit{A change of variable $z\mapsto w$ will be said to be admissible if it has the form
$$
w=z+P(x_1,y_1,\dots,x_n,y_n)
$$ 
where $P$ is a polynomial of degree 2 with no constant term.}

It is clear that the differential equation satisfied by $w$ has again the form
$$
\dot{w}=dw+l'(x_1,y_1,\dots,x_n,y_n)+Q'(x_1,y_1,\dots,x_n,y_n),
$$
where $l'$ is a linear form and  $Q'$ is a quadratic form in the $2n$ variables $x_1,y_1,\dots,x_n,y_n$.

\newtheorem{Formequadratiquebis}[champ linearsimplyconnex]{Theorem}
\begin{Formequadratiquebis}\label{Formequadratiquebis}
It is assumed that $Z\notin \LL\A\{B_1,\dots,B_m\}$ and $d\neq 0$.

If for some admissible change of variable the quadratic form $Q'$ is non negative and $\ker(l')\subset \ker(Q')$, then the system is not controllable.
\end{Formequadratiquebis}

\demo
Just repeat the arguments of the proof of Theorem \ref{Formequadratique}.

\hfill $\Box$

The previous paper \cite{DJ14} dealt with the $3$-dimensional case $\Hh^1$, and the proof of Theorem \ref{theohei3} (quoted in Section \ref{ContHun}) in the regular case $d\neq 0$ consists in proving the converse of Theorem \ref{Formequadratiquebis}. More accurately:

\textit{In $\Hh^1$, under the rank condition, a regular system is controllable if and only if for no admissible change of variable holds the condition: the quadratic form $Q'$ is non negative, and $\ker(l')\subset \ker(Q')$.}

This result was proved by considering the worst (in some sense) change of variable.

\vskip 0.2cm

It is consequently natural to make the following conjecture:

\vskip 0.2cm

\noindent \textbf{Conjecture}
\textit{Let $(\Sigma)$ be a regular system (that is $d\neq 0$) in $\Hn$ and assume that $Z\notin \LL\A\{B_1,\dots,B_m\}$.}

\textit{Then $(\Sigma)$ is controllable if and only if for no admissible change of variable holds the condition: the quadratic form $Q'$ is non negative and $\ker(l')\subset \ker(Q')$.}


\subsection{Sufficient controllability conditions}

The first proposition relates the rank condition on a quotient group to the rank condition of $(\Sigma)$.

\newtheorem{P16}[AlgebreLie]{Proposition}
\begin{P16}\label{P16}
Let $\g$ be an ideal of $\h^n$ invariant by $D$, and let $G$ be the subgroup generated by $\g$. Then $G$ is a closed Lie subgroup of $\Hn$ and the quotient $\mathbb{H}^{n}/G$ is an Abelian simply connected Lie group.

The induced system on $\mathbb{H}^{n}/G$ satisfies the rank condition as soon as $(\Sigma)$ does. It is in that case controllable in exact time $T$ for any $T>0$.
\end{P16}
\demo

The first assertion is proved in any textbook on Lie groups (see for instance \cite{Bourbaki2}).

Let us then assume that $(\Sigma)$ satisfies the rank condition, in other words that:
$$
\LL\A\{D^kB_j;\ \ k=0,\dots,2n,\ \ j=1,\dots,m\}=\hn.
$$

Let $\Pi$ stand for the projection of $\h^n$ onto $\h^n/\g$. Since it is a Lie algebra morphism, and since the algebra $\h^n/\g$ is Abelian, one has:
$$
\begin{array}{l}
\mbox{Vect}\{\Pi D^kB_j;\ \ k=0,\dots,2n,\ \ j=1,\dots,m\}\\
\qquad \qquad \qquad = \LL\A\{\Pi D^kB_j;\ \ k=0,\dots,2n,\ \ j=1,\dots,m\}\\
\qquad \qquad \qquad = \Pi\LL\A\{ D^kB_j;\ \ k=0,\dots,2n,\ \ j=1,\dots,m\}\\
\qquad \qquad \qquad =\Pi \h^n=\h^n/\g.
\end{array}
$$
This proves that the rank condition is satisfied on the quotient $\mathbb{H}^{n}/G$ as soon as it holds on $\Hn$. Since this quotient is simply connected and Abelian, the induced system is linear in the classical sense. It is consequently controllable in time $T$ for all $T>0$.

\hfill $\Box$

\vskip 0.2cm

\newtheorem{CasSingulier}[champ linearsimplyconnex]{Theorem}
\begin{CasSingulier}\label{CasSingulier}
If the linear system satisfies the algebraic rank condition and is singular (that is $d$, defined by $DZ=dz$, vanishes) then it is controllable.
\end{CasSingulier}
\demo
Since $\Sigma$ satisfies the rank condition the system induced on $\mathbb{H}^{n}/ \mathcal{Z}(\mathbb{H}^{n})$ is controllable. 

Since the algebraic rank condition holds the linearization at $e$ of the extended system is controllable (see Section \ref{ContLocCondAlg}) and the sets $\overline{\A}$ and $\overline{\A^-}$ contain neighbourhoods of $e$. But the system is real analytic and the interior of $\A$ (resp. of $\A^-$) is equal to the interior of $\overline{\A}$ (resp. of $\overline{\A^-}$).
Consequently $\A \cap \mathcal{Z}(\mathbb{H}^{n}) $ and $\A^{-} \cap \mathcal{Z}(\mathbb{H}^{n}) $
are neighbourhoods of $e$ in $\mathcal{Z}(\mathbb{H}^{n})$. Since the system is singular the linear vector field $\X$ vanishes on $\mathcal{Z}(\mathbb{H}^{n})$ and, according to Proposition \ref{propgen2}, $\mathcal{Z}(\mathbb{H}^{n})$ is included in  $\A$  and  $\A^{-}$. By Proposition \ref{propgen1} the system is controllable.

\hfill $\Box $

\vskip 0.5cm

\newtheorem{ZdansLA}[champ linearsimplyconnex]{Theorem}
\begin{ZdansLA}\label{ZdansLA}
If the rank condition is satisfied and if the invariant vector field $Z$ belongs to the Lie algebra generated by $ B_{1},\cdots, B_{m}$, then the system is controllable in exact time $T$ for all $ T > 0 $.
\end{ZdansLA}
\demo

As in the previous proof the system induced on $\mathbb{H}^{n}/ \mathcal{Z}(\mathbb{H}^{n})$ is controllable. 

If $ Z \in \LL\A \lbrace B_{1},\cdots, B_{m}\rbrace $, then $vZ$ belongs to the strong Lie saturate for all $v\in\R$ and we can consider the extended system
$$
(\Sigma_{E})\qquad \dot{g} = \X(g) + \sum_{j = 1}^{m}u_{j}B_{j}(g) + vZ.
$$
For $u_{j}=x_{i}=y_{i}=0$,  $i\in \left\lbrace 1,...n \right\rbrace $ and $j\in \left\lbrace 1,...m \right\rbrace $ the system reduced to $\mathcal{Z}(\mathbb{H}^{n})$ writes:
$\dot{z}=dz + v$. It is obviously controllable, so that  $\mathcal{Z}(\mathbb{H}^{n})$ is included in $\overline{\A} $ et $\overline{\A^{-}}$.

According to Proposition \ref{propgen1} the system is controllable on $\mathbb{H}^{n}$.\\

Let now $T > 0$ and $g=(x_{1}, y_{1}, x_{2}, y_{2},.....x_{n}, y_{n},z_{1})$ a point of $\mathbb{H}^{n}$. There exist admissible  controls $ u_{j}$, $j\in\lbrace 1, \cdots, m \rbrace $  that steer the identity $e$ to the class of $g$ in the quotient  $\mathbb{H}^{n}/ \mathcal{Z}(\mathbb{H}^{n})$ in time $\frac{T}{2}$. In $\mathbb{H}^{n}$, these controls steer the point $(0, 0, 0, ,.....0, 0, z_{2})$ to  $g=(x_{1}, y_{1}, x_{2}, y_{2},.....x_{n}, y_{n},z_{1})$ for some real number $z_2$. But in $\mathbb{H}^{n}$, there exist controls that steer the identity to $(0, 0, 0, ,.....0, 0, z_{2})$ in time $\frac{T}{2}$. The point $g$ can consequently be reached from $e$ in time $T$. Similarly the point  $e$ can be reached from $g$ in time  $T$, which proves that $(\Sigma)$
is controllable in time $T$ for all $T>0$.

\hfill $\Box$

Thanks to Theorem \ref{ZdansLA} we can now assert that the system is controllable as soon as $m\geq n+1$.
\newtheorem{C1}{Corollary}
\begin{C1}\label{C1}
The controlled vectors $B_1,\dots,B_m$ are assumed to be linearly independant and the rank condition to hold.

If $ m \geq n+1 $, then $(\Sigma)$ is controllable in exact time $T$ for all $ T > 0$. 
\end{C1}

\demo
According to Theorem \ref{ZdansLA} it is sufficient to prove that $Z$ belongs to the Lie algebra generated by $B_1,\dots,B_m$.

If $Z$ belongs to the linear space generated by $B_1,\dots,B_m$ the proof is finished. If not at least one of the brackets $[B_i,B_j]$ does not vanish. Indeed, if all these brackets would vanish we could choose a symplectic basis such that $X_i=B_i$ for $i=1,\dots,n$. Let 
$$
B_{n+1} = \sum_{i =1}^{n} \alpha_{i}B_{i} + \sum_{i =1}^{n} \beta_{i}Y_{i} + \gamma Z.
$$
Since $Z$ does not belong to the linear space generated by $B_1,\dots,B_m$, and since these vector fields are independant, at least one of the $\beta_{i}$ does not vanish, which implies that the brackets $[B_i,B_{n+1}]$ are not all zero, a contradiction.

\hfill $\Box$


\section{Decoupling}\label{Decoupling}
We consider again the system
$$
(\Sigma) \qquad \dot{g}=\X(g)+\sum_{j=1}^m u_j B_j(g)
$$
where $m\geq 1$ and the controlled vectors $B_1,\dots,B_m$ are linearly independant.

We have seen that $(\Sigma)$ is controllable as soon as one of the brackets $[B_i,B_j]$ is non zero (under the rank condition). In this section we investigate the case where
$$
[B_i,B_j]=0 \qquad \qquad \forall\ 1\leq i,j \leq m.
$$
According to Corollary \ref{C1} we can assume $m\leq n$ and according to Proposition \ref{Symplec Vectors} there exists a symplectic basis $(X_1,Y_1,\dots,X_n,Y_n,Z)$ such that $X_i=B_i$ for $i=1,\dots,m$.

In what follows we call $j^{th}$ \textit{cell} the bidimensional subspace of $\g$ generated by $B_j=X_j$ and $Y_j$, for $j\leq m$. Notice that the linear space generated by $(B_j,Y_j,Z)$ is an ideal of $\g$.

\newtheorem{decouplage}[syslin]{Definition}
\begin{decouplage} \label{decouplage}
The $j^{th}$ cell is said to be  \textit{decoupled} if the linear space generated by $(B_j,Y_j,Z)$ is invariant by $D$.
\end{decouplage}

This means that the matrix of $D$ has the following property: in the columms of $DX_j$ and $DY_j$ all the coefficients are zero excepted the ones corresponding to $X_j$, $Y_j$, and $Z$. Thanks to the symmetry of the derivations of $\hn$  all the coefficients of the lines of $X_j$ and $Y_j$ are zero excepted the ones corresponding to $X_j$, $Y_j$, and the differential equations satisfied by $x_j$ and $y_j$ does not depend on the other coordinates.

Let $G_j$ be the Lie subgroup generated by $(B_j,Y_j,Z)$. If the $j^{th}$ cell $(B_j,Y_j)$ is decoupled then on the first hand the subgroup $G_j$ is stable by the flow of $\X$ and $(\Sigma)$ induces a system on the quotient $\Hn/G_j$. On the other hand the behaviour of the two first coordinates of $G_j$, that is $(x_j,y_j)$, does not depend on the behaviour of the other parts of $(\Sigma)$. The restriction of $(\Sigma)$ to $G_j$ will be denoted by $(\Sigma_j)$.

\newtheorem{cellule cond rang}{Lemma}
\begin{cellule cond rang} \label{cellule cond rang}
If $(\Sigma)$ satisfies the rank condition, and if the $j^{th}$ cell is decoupled, then the system $(\Sigma_j)$ on $G_j$ satisfies also the rank condition.
\end{cellule cond rang}
\demo
Let us assume that $(\Sigma_j)$ does not satisfy the rank condition. This is equivalent to the fact that $DB_j=\alpha B_j$ for some constant $\alpha$. Then $Y_j$ does not belong to $\mbox{Span}\{D^kB_j;\ k\geq 1\}$. But the $j^{th}$ cell is decoupled and $Y_j$ neither belongs to $\mbox{Span}\{D^kB_i;\ i\neq j,\ k\geq 1\}$. Therefore $Y_j$ does not belong to the system Lie algebra and $(\Sigma)$ does not satisfy the rank condition.

\hfill $\Box$

\vskip 0.2cm

\newtheorem{cellule controlable}[champ linearsimplyconnex]{Theorem}
\begin{cellule controlable} \label{cellule controlable}
The rank condition is assumed to hold on $\Hn$.

If the $j^{th}$ cell is decoupled, and if $(\Sigma_j)$ is controllable on $G_j$ then $(\Sigma)$ is controllable on $\Hn$.
\end{cellule controlable}

\demo
Since $G_j$ is invariant by the flow of $\X$, $(\Sigma)$ induces a system on $H^n/G_j$. The quotient $H^n/G_j$ is the Abelian simply connected Lie group $\R^{2n-2}$, the induced system is linear and satisfies the rank condition, hence it is controllable.

On the other hand the subgroup $G_i$ is included in $\A$ and $\A^-$. Indeed let $g\in G_i$. Since $(\Sigma_j)$ is controllable by assumption, there exists a control $t\mapsto u_j(t)$ on some time interval $[0,T]$ that steers $e$ to $g$ in $G_j$. But it is clear that the control defined by $u_i=0$ if $i\neq j$ steers $e$ to $g$ in $\Hn$, which proves that $G_j\subset \A$.

Similarly $G_j\subset \A^-$ and according to Proposition \ref{propgen1} the system $(\Sigma)$ is controllable.

\hfill $\Box$


\subsection{Decoupled systems in $\Hh^2$}

In order to investigate in more details the dimension $5$, i.e. decoupled systems in $\Hh^2$, we first exhibit normal forms for these systems
\newtheorem{FormeNormaleDecouple}[cellule cond rang]{Lemma}
\begin{FormeNormaleDecouple}\label{FormeNormaleDecouple}
Let $(\Sigma)$ be a two inputs decoupled system in $\mathbb{H}^2$. If it satisfies the rank condition, there exist a symplectic basis $\{B_1,Y_1,B_2,Y_2,Z\}$ such that the matrix of the derivation $D$ be: 
$$
D=\begin{pmatrix}
0 & b & 0 & 0 & 0 \\
1 & d & 0 & 0 & 0 \\
0 & 0 & 0 & b' & 0 \\
0 & 0 & c' & d & 0 \\
0 & f & 0 & f' & d \\
\end{pmatrix}
\qquad \mbox{with}\quad c'\neq 0.
$$
\end{FormeNormaleDecouple}

\demo
The system being decoupled, the matrix of $D$ has the following form, in any symplectic basis such that $X_1=B_1$ and $X_2=B_2$:
$$
D=\begin{pmatrix}
a & b & 0 & 0 & 0 \\
c & d-a & 0 & 0 & 0 \\
0 & 0 & a' & b' & 0 \\
0 & 0 & c' & d-a' & 0 \\
e & f & e' & f' & d \\
\end{pmatrix}.
$$
The rank condition forces $c\neq 0$ and $c'\neq 0$.

First of all we can replace $Y_1$ by $B_1$, but in order that the basis remains symplectic we should
\begin{enumerate}
\item replace $Z$ by $cZ$, because $[B_1,DB_1]=cZ$;
\item replace $Y_2$ by $cY_2$ in order that $[B_2,cY_2]=cZ$.
\end{enumerate}
In the new basis we get:
$$
D=\begin{pmatrix}
0 & b & 0 & 0 & 0 \\
1 & d & 0 & 0 & 0 \\
0 & 0 & a' & b' & 0 \\
0 & 0 & c' & d-a' & 0 \\
0 & f & e' & f' & d \\
\end{pmatrix}
$$
(the parameters need not be the same than in the previous basis).

We have now $[B_2,DB_2]=c'Z$ but we can no longer modify $Z$. Consequently we cannot replace $Y_2$ by $DB_2$ but only by $\displaystyle \frac{1}{c'} DB_2$, which gives the matrix stated in the Lemma.

\hfill $\Box$

\newtheorem{dimCinq}[champ linearsimplyconnex]{Theorem}
\begin{dimCinq}\label{dimCinq}
Let $(\Sigma)$ be a two inputs decoupled system in $\mathbb{H}^2$ in the form of Lemma \ref{FormeNormaleDecouple}. It is assumed to satisfy the rank condition.

If none of the subsystems $(\Sigma_1)$ and $(\Sigma_2)$ is controllable, then $(\Sigma)$ is controllable if and only if $c'$ is negative.
\end{dimCinq}

\demo

The system being in the form of Lemma \ref{FormeNormaleDecouple} it writes in the associated coordinates:
$$
(\Sigma)=
\left\{
\begin{array}{ll}
\dot{x_1} & =by_1+u_1\\
\dot{y_1} & =x_1+dy_1\\
\dot{x_2} & =b'y_2+u_2\\
\dot{y_2} & =c'x_2+dy_2\\
\dot{z} & =fy_1+f'y_2+dz+\frac{1}{2}x_1^2+\frac{1}{2}by_1^2+\frac{1}{2}c'x_2^2+\frac{1}{2}b'y_2^2
\end{array}
\right.
$$
Let us first assume that $c'>0$ and let us show that $(\Sigma)$ is not controllable.

Since the cells are not controllable, and according to the results of \cite{DJ14}, we have:
$$
\displaystyle b\geq -\frac{d^2}{4}, \qquad b'\geq -\frac{d^2}{4c'},\ \  \mbox{ and if }\ \ d=0\ \ \mbox{ then } f=f'=0.
$$
The variable $w=z+\frac{d}{4}y_1^2+\frac{d}{4c'}$ satisfies the differential equation:
$$
\dot{w}=fy_1+f'y_2+dw+\frac{1}{2}(x_1+\frac{d}{2}y_1)^2+\frac{1}{2}(b+\frac{d^2}{4})y_1^2+\frac{c'}{2}(x_2+\frac{d}{2c'}y_2)^2
+\frac{1}{2}(b'+\frac{d^2}{4c'})y_2^2
$$

If $\displaystyle b> -\frac{d^2}{4}$ and $\displaystyle b'> -\frac{d^2}{4c'}$, the quadratic part of this equation is positive definite and according to Theorem \ref{Formequadratiquebis} the system cannot be controllable.

In case of equality, for instance if $\displaystyle b= -\frac{d^2}{4}$, and if $d\neq 0$, the second change of variable  $\displaystyle w\longmapsto w+2\frac{f}{d}y_1$ leads to the same conclusion (see \cite{DJ14} or \cite{Dath15} for more details).

If $d=0$ the non controllability of the cells implies $f=f'=0$ hence:
$$
\dot{w}=\frac{1}{2}x_1^2+\frac{1}{2}by_1^2+\frac{c'}{2}x_2^2
+\frac{1}{2}b'y_1^2.
$$
This quadratic form is positive for $b\geq 0$ and $b'\geq 0$ and $(\Sigma)$ is not controllable.

\vskip 0.2cm

Let us now show that $c'<0$ implies controllability.

For $i=1,2$ the vector field $vB_i$ belongs to $\LS(\Sigma)$ for all $v\in \R$, hence the vector field $\exp(v\ad(B_i))\X$ belongs as well to $\LS(\Sigma)$ (see \cite{Jurdjevic97}). But
$$
\begin{array}{ll}
\exp(v\ad(B_i))\X&=\X+(v\ad(B_i))\X + \frac{1}{2}(v\ad(B_i))^{2} \X \\
&= \X + v DB_i +\frac{v^{2}}{2}[B_i,DB_i].
\end{array}
$$
For $i=1$, $\X + v DB_1 +\frac{v^{2}}{2}[B_1,DB_1]=\X + v DB_1 +\frac{v^{2}}{2}Z$, hence the vector field $Z$ belongs to the system Lie saturate (dividing by $v^2$, we get $Z$ when $v$ tends to $+\infty$).

For $i=2$, $\X + v DB_2 +\frac{v^{2}}{2}[B_2,DB_2]=\X + v DB_2 +\frac{v^{2}}{2}c'Z$, and the vector $c'Z$ belongs to the Lie saturate.

Since $c'<0$ we finally obtain that $\pm Z\in\LL\Sg$, hence that the center of $\Hh^2$ is included in $\overline{\A}$ and  $\overline{\A^-}$, and according to Proposition \ref{propgen1} that $(\Sigma)$ is controllable.

\hfill $\Box$


\subsection{Extension to $\Hn$}

Consider in $\Hn$ a system with $m$ inputs and that satisfies the rank condition. Assume that $m_1$ cells, with $2\leq m_1\leq m$, are decoupled. We can assume that the decoupled cells are the $m_1$ first ones and we can choose a symplectic basis such that $X_i=B_i$ for $i=1,\dots,m_1$ and such that the restriction of $D$ to $(B_i,Y_i,Z)$ is:
$$
\begin{pmatrix}
0 & b_i & 0 \\
c_i & d & 0 \\
0 & f_i & d 
\end{pmatrix}
\qquad \mbox{ with }\ \ c_i\neq 0
$$
We can moreover assume that $c_1=1$.

According to the previous section the vector field $+Z$ belongs to the Lie saturate because $c_1=1$. If one of the $c_i$'s ($i=2,\dots,m_1$) is negative the vector field $-Z$ also belongs to the Lie saturate and the system is controllable. If all the $c_i$'s are positive, we cannot conclude to non controllability because we do not know the behaviour of the system with respect to the controls $u_i$ for $i=m_1+1,\dots, m$. Notice that if $m_1<n$ then $m_1<m$. Indeed if all the cells are decoupled, the rank condition cannot be satisfied if $m<n$. This remark leads to what follows.

Consider in $\Hn$ a system with $m=n$ inputs whose all cells are decoupled and that satisfies the rank condition. Such a system will be refered to as a \textit{completely decoupled system}.

With the same notations as above we obtain at once that in the case where all the $c_i$'s are positive the same kind of change of variables than in the proof of Theorem $\ref{dimCinq}$ provides a positive definite quadratic form that will be an obstruction to controllability.

We can therefore state a generalization of Theorem \ref{dimCinq}.

\newtheorem{dimPlusqueCinq}[champ linearsimplyconnex]{Theorem}
\begin{dimPlusqueCinq}\label{dimPlusqueCinq}

Let $(\Sigma)$ be a system with $m$ inputs, with $m_1$ decoupled cells, and that satisfies the rank condition.

If none of the subsystems $(\Sigma_i)$  is controllable, but if one of the $c_i$'s is negative then $(\Sigma)$ is controllable.

If the system is completely decoupled and if none of the subsystems $(\Sigma_i)$ is controllable, then $(\Sigma)$ is controllable if and only if one of the $c_i$'s is negative.
\end{dimPlusqueCinq}


\section{Conclusion}

In the previous paper \cite{DJ14} we had shown that under the rank condition a regular system on $\Hh^1$ is controllable if and only if the eigenvalues of the system induced on $\Hh^1/\Z(\Hh^1)$ are not real (see Theorem \ref{theohei3} and Corollary \ref{coroheis} in Section \ref{ContHun}).

Thanks to the results of Section \ref{Decoupling} about decoupling we know that this condition is no longer necessary on the generalized Heisenberg group $\Hn$ if $n\geq 2$. Indeed consider in $\Hn$ a regular system with two decoupled cells. If  their eigenvalues are real these cells are non controllable. However if one of the coefficients $c_i$ that appear in normal form is positive and the other one negative, then $(\Sigma)$ is controllable despite the existence of real eigenvalues in the quotient.


\end{document}